\documentclass[twoside]{amsart}
\usepackage{latexsym}
\usepackage{amssymb,amsmath,amsopn}
\usepackage[dvips]{graphicx}   
\usepackage{color,epsfig}      

\newtheorem*{thm}{Theorem}
\newtheorem*{lem}{Lemma}
\theoremstyle{definition}

\newtheorem*{prob}{Problem}

\renewcommand{\Re}{\mathbb R}
\newcommand{\Eu}{\mathbb E}
\newcommand{\M}{\mathbb{M}}

\newcommand{\HH}{\mathbb{H}}

\DeclareMathOperator{\inter}{int}
\DeclareMathOperator{\bd}{bd}
\DeclareMathOperator{\dist}{dist}
\DeclareMathOperator{\conv}{conv}

\DeclareMathOperator{\perim}{perim}
\DeclareMathOperator{\area}{area}
\DeclareMathOperator{\diam}{diam}

\begin{document}
\title[On the perimeters of simple polygons]{On the perimeters of simple polygons contained in a plane convex body}

\author[Z. L\'angi]{Zsolt L\'angi}
                                                                               
\address{Zsolt L\'angi, Dept. of Geometry, Budapest University of Technology,
Budapest, Egry J\'ozsef u. 1., Hungary, 1111}
\email{zlangi@math.bme.hu}

\thanks{The author was supported by the J\'anos Bolyai Research Scholarship of the Hungarian Academy of Sciences.}
                                                                        

\begin{abstract}
A simple $n$-gon is a polygon with $n$ edges such that each vertex belongs to exactly two edges and
every other point belongs to at most one edge.
Brass, Moser and Pach \cite[Problem 3, p. 437]{BMP05} asked the following question: For $n \geq 5$ odd,
what is the maximum perimeter of a simple $n$-gon contained in a Euclidean unit disk?
In 2009, Audet, Hansen and Messine \cite{AHM09} answered this question, and showed that the
supremum is the perimeter of an isosceles triangle inscribed in the disk, with an edge of multiplicity $n-2$.
In \cite{L11}, L\'angi generalized their result for polygons contained in a hyperbolic disk.
In this note we find the supremum of the perimeters of simple $n$-gons contained in an arbitrary convex body
in the Euclidean or in the hyperbolic plane.
\end{abstract}
\maketitle

\section{Introduction}

A question in the spirit of isoperimetric problems about simple polygons
was asked by Brass, Moser and Pach (cf. \cite[Problem 3, p. 437]{BMP05}).

\begin{prob}[Brass, Moser and Pach, 2005]
For $n \geq 5$ odd, what is the maximum perimeter of a simple $n$-gon contained in a Euclidean unit disk?
\end{prob}

The authors of \cite{BMP05} remarked that for $n$ even, the supremum of the perimeters is
the trivial upper bound $2n$, as it can be approached by simple $n$-gons in which the vertices alternate
between some small neighborhoods of two antipodal points of the disk.
This argument cannot be applied if $n$ is odd.
In 2009, Audet, Hansen and Messine \cite{AHM09} showed that for $n$ odd, the supremum is attained
by the perimeter of an isosceles triangle inscribed in the disk, with an edge of multiplicity $n-2$.
The author of \cite{L11} gave a shorter proof of the same statement and proved that
for hyperbolic disks of any radius, the supremum is attained by the perimeter of an $n$-gon of the same kind;
that is, by the perimeter of an isosceles triangle with a multiple edge inscribed in the disk.
He noted that for $n$ even and for any convex body $C$ in the Euclidean plane $\Eu^2$ or in the hyperbolic plane $\HH^2$,
the supremum of the perimeters of the simple $n$-gons contained in $C$ is the trivial bound $n \diam C$, where $\diam C$ is the diameter of $C$.
He asked whether it is true that, for $n$ odd, the supremum is the perimeter of a triangle with an edge
of multiplicity $n-2$, inscribed in $C$.

In this paper we answer this question. Our main result is the following.

\begin{thm}\label{thm:oval}
Let $n \geq 3$ be an odd integer, and let $C$ be a convex body in $\Eu^2$ or in $\HH^2$.
For every simple $n$-gon $P$ contained in $C$ there is a triangle, inscribed in $C$ and with side-lengths $\alpha \geq \beta \geq \gamma$, such that $\perim P \leq (n-2)\alpha + \beta + \gamma$.
\end{thm}

\begin{figure}[here]
\includegraphics[width=0.6\textwidth]{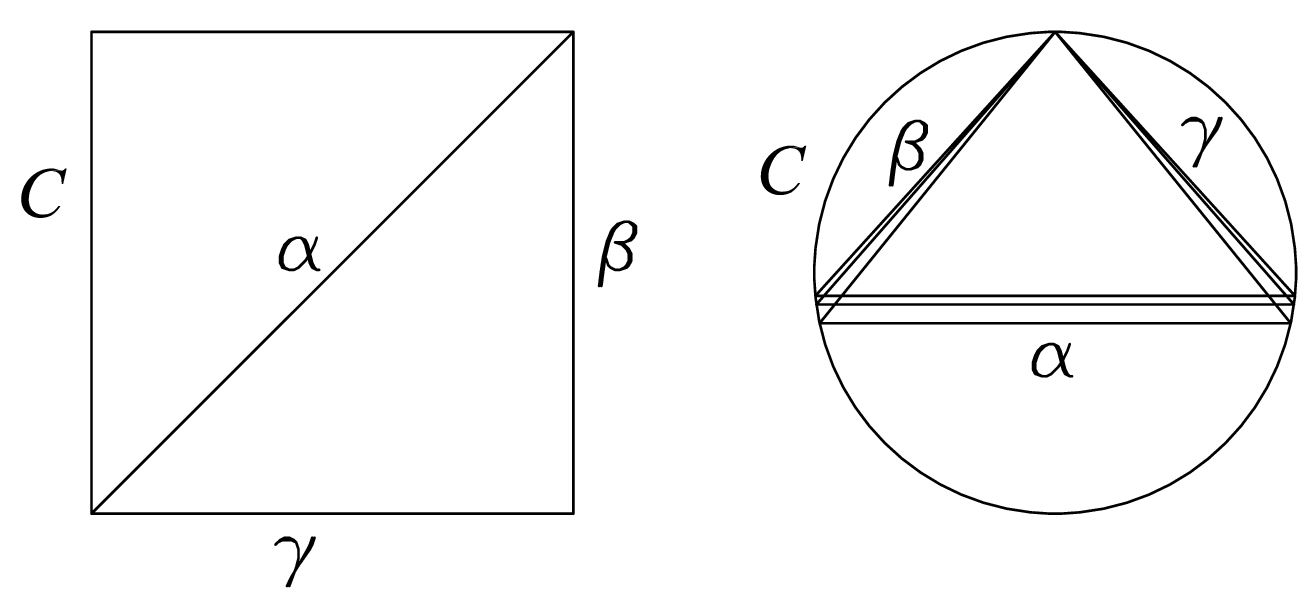}
\caption[]{}
\label{fig:example}
\end{figure}

Figure~\ref{fig:example} shows such triangles of maximum perimeter
for a square for any value of $n$, and for a unit disk for $n=5,7,9$. We note that the side-lengths of this triangle
are independent of $n$ if $C$ is a square, and that their values for a unit disk were determined in \cite{AHM09}.

In the proof we use the following notations.
Let $\M \in \{ \Eu^2, \HH^2\}$ and $x, y\in \M$.
The distance of $x$ and $y$ is denoted by $\dist(x,y)$.
The closed (respectively, open) segment with endpoints $x$ and $y$ is denoted by $[x,y]$ (respectively, $(x,y)$).
If $x \neq y$, $L(x,y)$ denotes the straight line passing through $x$ and $y$, and $R_x(x,y)$
denotes the closed ray in $L(x,y)$ emanating from $x$ and not containing $y$.

For any set $A \subset \M$, we use the standard notations
$\inter A$, $\bd A$, $\diam A$, $\perim A$, $\area A$ and $\conv A$ for the \emph{interior}, the \emph{boundary},
the \emph{diameter}, the \emph{perimeter}, the \emph{area}, or the \emph{convex hull} of $A$.
Points are denoted by small Latin letters, and sets of points by capital Latin letters.

In the proof we find a triangle $T$, \emph{contained} in $C$, with side-lengths $\alpha, \beta$
and $\gamma$ such that $\perim P \leq (n-2)\alpha + \beta + \gamma$,
as in this case we can move the vertices of $T$
to $\bd C$ in a way that no side-length of $T$ decreases.

\section{Proof of Theorem}\label{sec:Euclidean}

We prove the theorem only for the Euclidean plane, since for $\HH^2$ one only needs to apply the Euclidean argument in a slightly modified way.

Let us consider a Cartesian coordinate system.
If $z \in \Eu^2$ is an arbitrary point, by $z = (\mu, \nu)$ we mean that the $x$-coordinate of $z$ is $\mu$, and its $y$-coordinate is $\nu$.
Let $[a_0,a_1], [a_1,a_2], \ldots, [a_{n-1},a_n]$ denote the edges of $P$ such that $a_0 = a_n$, and let $a_i = (\omega_i,\theta_i)$ for $i=0,1,2,\ldots,n$.
Without loss of generality, we may assume that $[a_0,a_1]$ is a longest edge of $P$,
$a_0$ is the origin $(0,0)$, and that $a_1=(0,1)$.

For $i=0,1,\ldots,n$, let $\zeta_i=\theta_{i+1} - \theta_i$.
Note that $\zeta_0=\zeta_n =1$.
As $n$ is odd, the sequence $\{ \zeta_i \}$ consists of an even number of elements.
Thus, it has two consecutive elements, say $\zeta_{j-1}$ and $\zeta_j$, that are both nonnegative or nonpositive.
From this, we have that $\theta_{j-1} \leq \theta_j \leq \theta_{j+1}$, or that
$\theta_{j-1} \geq \theta_j \geq \theta_{j+1}$, respectively.
For simplicity, we denote $a_0$, $a_1$, $a_{j-1}$, $a_j$ and $a_{j+1}$ by $p=(0,0)$, $q=(0,1)$, $a=(\omega_a,\theta_a)$, $b=(\omega_b,\theta_b)$ and
$c=(\omega_c,\theta_c)$, respectively, and set $p_a=(0,\theta_a)$ and $p_c = (0,\theta_c)$.
We remark that this argument is used both in \cite{AHM09} and in \cite{L11}.

During the proof, we may, without loss of generality, assume that
\begin{enumerate}
\item $\theta_a \leq \theta_b \leq \theta_c$,
\item $a$ is not farther from the bisector of $[p,q]$ than $c$ (or in other words, $ \theta_a + \theta_c \geq 1$),
\item at least one of $\omega_a$, $\omega_b$ and $\omega_c$ is positive.
\end{enumerate}

A possible approach to prove the Theorem is to ignore all the edges of $P$ but $[p,q]$, $[a,b]$ and $[a,c]$,
and to show the existence of points $a',b',c' \in C$
that satisfy $\dist(p,q)\leq \dist(a',c')$ and $\dist(a,b)+\dist(a,c) \leq \dist(a',b')+\dist(b',c')$, from which the assertion would readily follow.
This was done in \cite{L11} for a Euclidean unit disk.
Unfortunately, this property does not hold for every plane convex body, as the following example shows.

\begin{figure}[!h]
\begin{minipage}[t]{0.3\columnwidth}%
\vspace{0pt}
    \centering
    \includegraphics[width=0.48\textwidth]{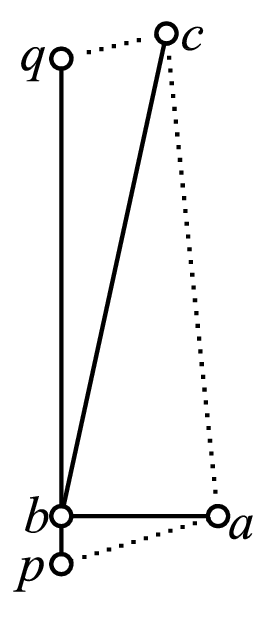}
    \caption[]{}
    \label{fig:pelda}
\end{minipage}
\hfill%
\begin{minipage}[t]{0.6\columnwidth}%
\vspace{7pt}

$p=(0,0)$, \quad $q=(0,1)$,\\
$a=(0.31,0.095)$, \quad $b=(0,0.095)$,\\
$c=(0.208,1.05)$,
$C = \conv \{ p,q,a,b,c\}$,
\bigskip

$\dist(a,b) = 0.3100 \ldots$,\\
$\dist(b,c) =  0.9773\ldots$,\\
$\dist(a,c) =  0.9604\ldots$,

\bigskip
$\dist(a,b)+\dist(b,c) = 1.2873\ldots$,\\
$\dist(p,c)+\dist(c,q) = 1.2843\ldots$,\\
$\dist(p,a)+\dist(a,q) = 1.2808\ldots$,\\
$\dist(p,a)+\dist(a,c) = 1.2846\ldots$.

\end{minipage}
\end{figure}

In the proof, we first show in the Lemma that the property described in the third paragraph of this section fails
only for fairly well-determined configurations, and then, after the Lemma, we prove the Theorem for these configurations in a slightly different way.
To prove the Lemma and the Theorem, we need geometric observations different to those used for a Euclidean unit disk.

For simplicity, if $a', b', c'$ satisfy $\dist(p,q)\leq \dist(a',c')$ and $\dist(a,b)+\dist(a,c) \leq \dist(a',b')+\dist(b',c')$,
we say that \emph{$a', b'$ and $c'$ satisfy Property (*)}.

\begin{lem}
If there are no points $a',b',c' \in \conv\{ p,q,a,b,c\}$ satisfying Property (*), then the following hold.
\begin{enumerate}
\item[(a)] $\dist(a,c) < 1$.
\item[(b)] $a, b$ and $c$ are in the same closed half-plane bounded by $L(p,q)$.
\item[(b)] $\theta_c > 1$ and $0 < \theta_a < \frac{1}{2}$. 
\item[(d)] $b \in \conv\{ p_a,p_c,a,c\}$.
\item[(e)] $\dist(a,b) + \dist(b,c) \leq \dist(a,p_a) + \dist(p_a,c)$.
\end{enumerate}
\end{lem}

\begin{proof}
We prove the Lemma, by contradiction, in six steps.
In the proof, each of the five conditions but (c) is proved in one step.
The remaining condition, (c), is proved in Steps 3 and 6.
Within each step, we already use the conditions proved in the previous steps.

\noindent
\emph{Step 1}.\\
If $\dist(a,c) \geq 1$, then $a, b$ and $c$ clearly satisfy Property (*); a contradiction.

\noindent
\emph{Step 2}.\\
Suppose, for contradiction, that $a$, $b$ and $c$ are not in the same closed half-plane bounded by $L(p,q)$.

First, we examine the case that $[a,b] \cap R_q(p,q) \neq \emptyset$.
Then $\dist(a,b) \leq 1 \leq \dist(b,p)$ and $1 \leq \dist(c,p)$
(cf. Figure~\ref{fig:figure1}), and thus, $p$, $b$ and $c$ satisfy Property (*).
If $[b,c] \cap R_p[p,q] \neq \emptyset$, we may apply a similar argument.
 
Next, we consider the case that $[a,b] \cap R_p(p,q) \neq \emptyset$, which yields that $\theta_a \leq 0$.
If $\theta_b \leq 0$, then we may apply the argument in the previous paragraph, and thus we have $0 < \theta_b$.
From this and from $\dist(a,c) < 1$, we readily obtain that $0 < \theta_b \leq \theta_c < 1$.
Let $L$ denote the bisector of the segment $[c,q]$.
Since $\dist(a,c) < 1 \leq \dist(a,q)$, we have that $L$ separates $q$ from $a$ and $c$. 
Observe also that $L$ separates $b$ and $q$ from $c$, as otherwise $\dist(b,q) \geq \dist(b,c)$, and $a$, $b$
and $q$ satisfy Property (*); a contradiction.
Thus, $[a,b] \cap L \neq \emptyset \neq [b,c] \cap L$, which implies that $b \in \conv \{ a, p_c, c \}$ (cf. Figure~\ref{fig:extra1}).
From this, we obtain that $\dist(a,b) + \dist(b,c) \leq \dist(a,p_c) + \dist(p_c,c) \leq \dist(a,c) + \dist(c,q)$,
and then $a, c$ and $q$ satisfy Property (*); a contradiction.
The case $[b,c] \cap R_q(p,q) \neq \emptyset$ follows by a similar argument.

\begin{figure}[!h]
\begin{minipage}[b]{0.49\columnwidth}%
    \centering
    \includegraphics[width=0.65\textwidth]{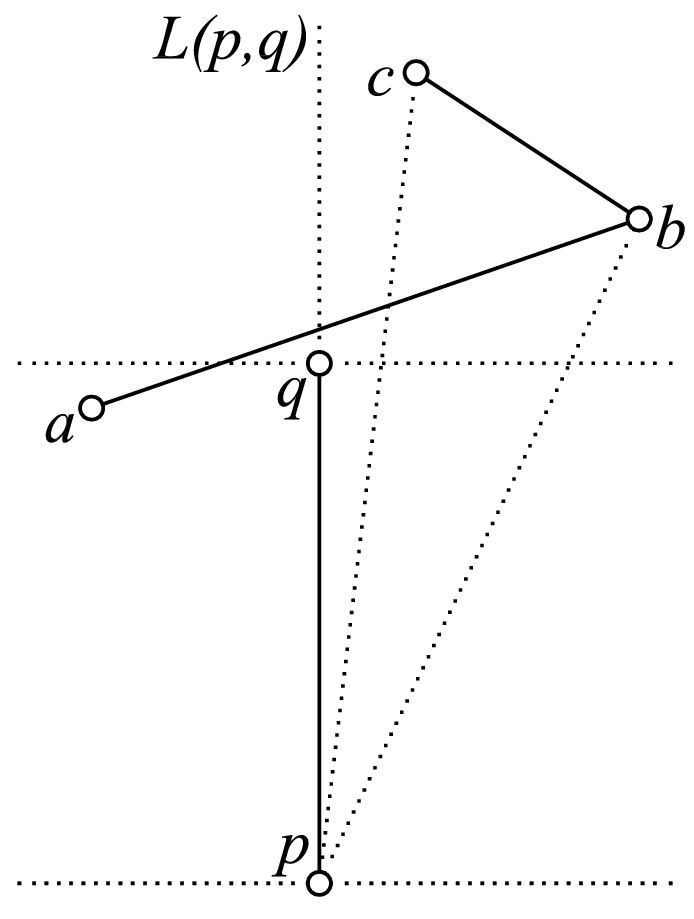}
		\caption[]{}
		\label{fig:figure1}
\end{minipage}
\hfill%
\begin{minipage}[b]{0.49\columnwidth}%
		\includegraphics[width=0.78\textwidth]{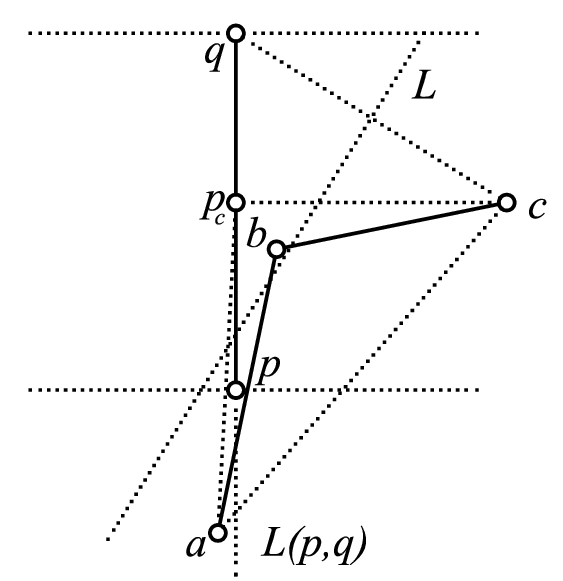}
		\caption[]{}
		\label{fig:extra1}
\end{minipage}
\end{figure}

In the rest of the proof we may assume that $a$, $b$ and $c$ are in the same closed half plane bounded by $L(p,q)$, which,
according to our assumptions, yields $0 \leq \omega_a$, $0 \leq \omega_b$ and $0 \leq \omega_c$.

\noindent
\emph{Step 3}.\\
Now we show that $\theta_c > 1$ and $0 < \theta_a < 1$.
First, observe that at least one of $0 < \theta_a$ and $\theta_c < 1$ holds, as otherwise $\dist(a,c) \geq 1$; a contradiction.
From this, as $\theta_a + \theta_c \geq 1$, it follows that $0 < \theta_a$.

Consider the case that $\theta_c \leq 1$, and let $\Omega = \max \{ \omega_a, \omega_b, \omega_c\}$.
If $\Omega = \omega_b$, then $\dist(a,b) \leq \dist(p_a,b) \leq \dist(b,p)$ and, similarly,
$\dist(b,c) \leq \dist(b,q)$, which yields that $p, b$ and $q$ satisfy Property (*).
If $\Omega = \omega_c$, then $\dist(a,b) + \dist(b,c) \leq \dist(p_a,c) + \dist(p_c,c)
\leq \dist(p,c)+\dist(c,q)$, and thus, $p, c$ and $q$ satisfy Property (*).
If $\Omega = \omega_a$, then the assertion follows by a similar argument.

Finally, if $\theta_a \geq 1$, then, by the argument
in Step 2, we have that $p$, $a$ and $c$ satisfy Property (*), which proves the last inequality.

\noindent
\emph{Step 4}.\\
Suppose for contradiction that $b \notin \conv \{ p_a,p_c,a,c\}$.
Then the three rays, emanating from $a$, that pass through $p$, $c$ and $b$ are in this clockwise order around $a$.
Let $L'$ denote the bisector of the segment $[p,a]$.
Note that as $\dist(a,c) < 1 \leq \dist(p,c)$, $L'$ separates $[a,c]$ from $p$.
Hence, it follows from $\theta_b \leq \theta_c$ that $L'$ separates $[a,b]$ from $p$
(cf. Figure~\ref{fig:figure2}).
Thus, $\dist(a,b) \leq \dist(p,b)$, and $p$, $b$ and $c$ satisfy Property (*); a contradiction.

\begin{figure}[!h]
\begin{minipage}[b]{0.49\columnwidth}%
    \centering
    \includegraphics[width=0.62\textwidth]{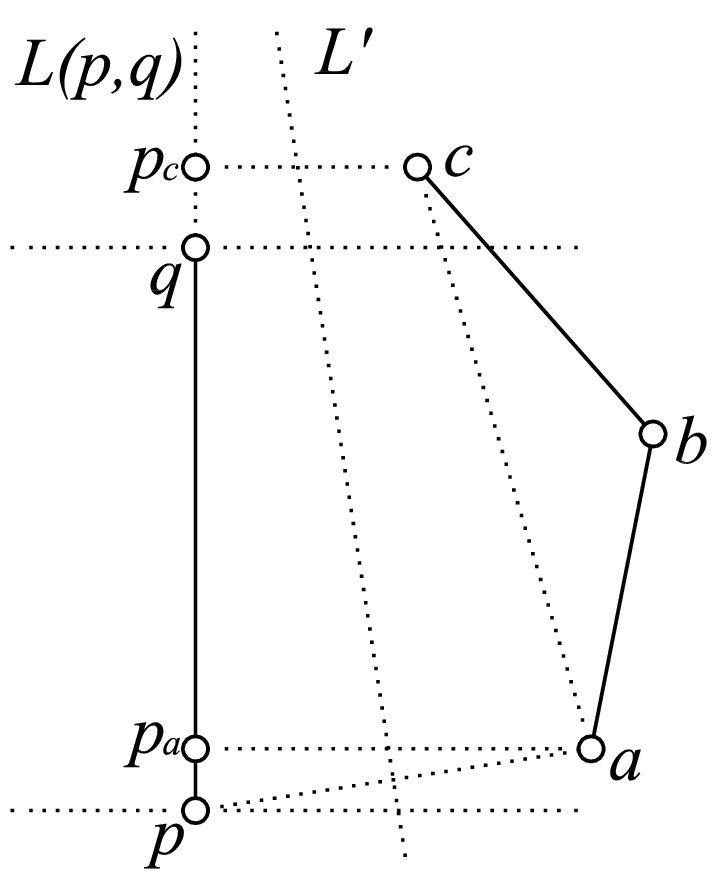}
		\caption[]{}
		\label{fig:figure2}
\end{minipage}
\hfill%
\begin{minipage}[b]{0.44\columnwidth}%
		\includegraphics[width=0.72\textwidth]{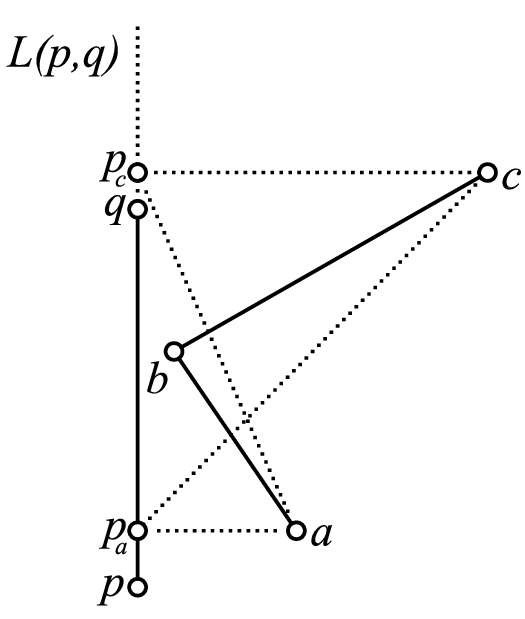}
		\caption[]{}
		\label{fig:extra2}
\end{minipage}
\end{figure}

\noindent
\emph{Step 5}.\\
We show that $\dist(a,b) + \dist(b,c) \leq \dist(a,p_a) + \dist(p_a,c)$.
If $b \notin \conv \{ p,q,a,c\}$ but $b \in \conv \{ p_a,p_c,a,c\}$, then $b \in \conv \{ q,p_c,c\}$, $\dist(b,c) \leq \dist(c,q)$
and $p, q$ and $c$ satisfy Property (*).
Thus, we have $b \in \conv \{ p,q,a,c\}$,  and this yields that
$\dist(a,b) + \dist(b,c) \leq \max \{ \dist(a,q) + \dist(q,c), \dist(a,p_a) + \dist(p_a,c) \}$.

Clearly, to prove our statement it suffices to deal with the case $\dist(a,p_a) + \dist(p_a,c) \leq \dist(a,q) + \dist(q,c)$.
Note that in this case $\dist(a,p_a) + \dist(p_a,c) \leq \dist(a,p_c) + \dist(p_c,c)$, which yields that
$\omega_a \leq \omega_c$ (cf. Figure~\ref{fig:extra2}).
Thus, the two legs of the right triangle $\conv \{ p,c,p_c\}$ are pairwise greater than or equal to the two legs of
$\conv \{ q,a,p_a\}$, from which we obtain that $\dist(q,a) \leq \dist(p,c)$,
and that $p, c$ and $q$ satisfy property (*).

\noindent
\emph{Step 6}.\\
Finally, we show that $\theta_a < \frac{1}{2}$.
For contradiction, assume that $\theta_a \geq \frac{1}{2}$.
This implies that $\dist(a,p) \geq \dist(a,q)$, and, by setting $u = \left( 0,\frac{1}{2} \right)$ and $v = \left( \omega_a, \frac{1}{2} \right)$,
that $\dist(a,b) + \dist(b,c) \leq \dist(a,p_a) + \dist(p_a,c) \leq \dist(u,v) + \dist(u,c)$.

Note that for any $x,y,z \in \Eu^2$, the function $\tau \mapsto \dist(x,y+\tau z)$ is a convex function
on $\Re$.
Thus,
\[
\dist(u,v) + \dist(u,c) \leq \frac{1}{2} \left( \dist(a,p) + \dist(a,q) \right)
 + \frac{1}{2} \left( \dist(p,c) + \dist(q,c) \right) \leq
\]
\[ 
\leq \max \{ \dist(a,p) + \dist(a,q), \dist(p,c) + \dist(c,q) \},
\]
and the assertion readily follows.
\end{proof}

Now we prove the Theorem.
Clearly, we may assume that the conditions from (a) to (e) hold for $p, q, a, b$ and $c$.
First, observe that $0 < \theta_a < \frac{1}{2}$ and $\theta_c > 1$ imply that $\theta_c - \theta_a > \frac{\theta_c}{2}$.
We distinguish two cases: $n = 5$ and $n \geq 7$.

First, let $n=5$. Then, under the conditions (a) to (e), we have that the 
remaining two edges of $P$ are $[p,a]$ and $[q,c]$.
Thus, $\perim P \leq 3\dist(p,c) + \dist(p,q) + \dist(q,c)$.

Second, we assume that $n \geq 7$.
If $\dist(a,c) \geq \dist(p_a,c)$, then $p, a$ and $c$ satisfy Property (*), from which the assertion readily follows.
Hence, we may assume that $\dist(a,c) \leq \dist(p_a,c)$, or in other words, that $\omega_c \geq \frac{\omega_a}{2}$.
Furthermore, since $n \geq 7$ and $1 = \dist(p,q) < \dist(p,c)$, it suffices to prove that
$5 + \dist(a,b) + \dist(b,c) \leq 5\dist(p,c) + \dist(a,p) + \dist(a,c)$; that is, that
\begin{equation}\label{eq:crucialeu}
5 + \dist(a,p_a) + \dist(p_a,c) \leq 5\dist(p,c) + \dist(a,p) + \dist(a,c).
\end{equation}
This is our aim for the remaining part of the proof.

Let $M$ and $N$ denote the left-hand side and the right-hand side of (\ref{eq:crucialeu}), respectively,
and let us regard them as functions of $c$.
Consider the vector $v=(1,0)$ and set $w = \left( \frac{\omega_a}{2}, \theta_c \right)$.
Note that $5 + \dist(a,p_a) + \dist(p_a,w) \leq 5\dist(p,w) + \dist(p,a) + \dist(a,w)$, which means that
(\ref{eq:crucialeu}) holds for $c=w$.
Thus, we need only to prove that in the direction of $v$, the derivative of $N$ is not smaller than that of $M$;
that is, using the standard notation from differential geometry, that $v(M) \leq  v(N)$.

\begin{figure}[here]
\includegraphics[width=0.34\textwidth]{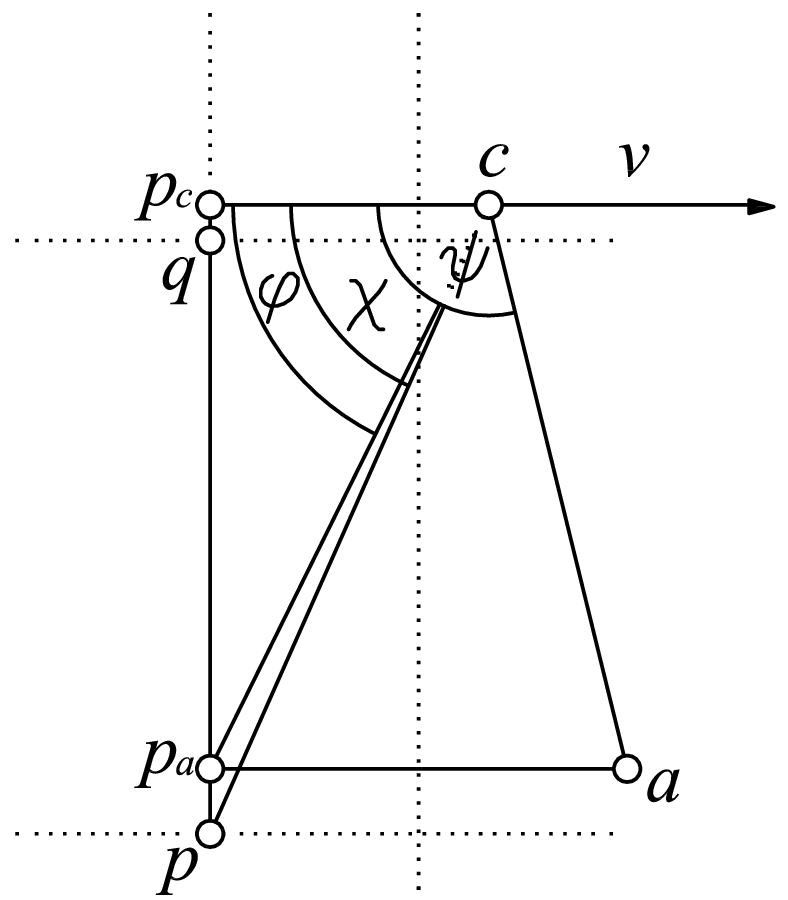}
\caption[]{}
\label{fig:figure3}
\end{figure}

Let $\phi, \chi$ and $\psi$ denote the internal angles at $c$ of $\conv \{ p_a,p_c, c\}$, $\conv \{ p,p_c,c\}$
and $\conv \{ a,c,p_a,p_c\}$, respectively (cf. Figure~\ref{fig:figure3}).
Observe that $ 0 < \phi \leq \pi - \psi < \pi$ and that $\cos \phi \geq - \cos \psi$.

Note that $v(M) = \cos \phi$, and $v(N) = 5 \cos \chi + \cos \psi \geq 5\cos \chi - \cos \phi$. 
We set $I = 5\cos \chi - 2 \cos \phi \leq v(N) - v(M)$.
Then an elementary calculation yields that
\[
I = \frac{5\omega_c}{\sqrt{\omega_c^2+\theta_c^2}}-\frac{2\omega_c}{\sqrt{\omega_c^2+(\theta_c-\theta_a)^2}} \geq
\frac{5\omega_c}{\sqrt{\omega_c^2+\theta_c^2}}-\frac{2\omega_c}{\sqrt{\omega_c^2+(\theta_c/2)^2}} =
\]
\[
= \frac{\omega_c\left( 21\omega_c^2+\frac{9}{4}\theta_c^2 \right)}{\sqrt{\omega_c^2+\theta_c^2}\sqrt{\omega_c^2+(\theta_c/2)^2}
\left( 5\sqrt{\omega_c^2+(\theta_c/2)^2}+2\sqrt{\omega_c^2+\theta_c^2} \right)} \geq 0,
\]
which finishes the proof of the theorem.

\end{document}